\documentclass[12pt]{article}
\usepackage{amsmath,amssymb}
\usepackage{enumerate}
\usepackage{a4wide}

\newcommand{\N}{{\mathbb N}}
\newcommand{\R}{{\mathbb R}}
\newcommand{\Q}{{\mathbb Q}}
\newcommand{\Z}{{\mathbb Z}}

\newcommand{\B}{\mathcal B}

\newcommand{\A}{\mathcal A}
\newcommand{\M}{\mathcal M}

\newcommand{\cL}{{\mathcal L}}
\renewcommand{\L}{{\mathcal L}}

\newtheorem{lemma}{Lemma}[section]
\newtheorem{proposition}[lemma]{Proposition}

\newtheorem{theorem}{Theorem}
\newcommand{\proof}{{\em Proof.}}
\newcommand{\cqfd}{\hfill$\Box$}

\newcommand{\esp}{{\mathbb E}}
\newcommand{\ind}{{\mathchoice{\mathrm {1\mskip-4.1mu l}} {\mathrm{ 1\mskip-4.1mu l}} {\mathrm {1\mskip-4.6mu l}} {\mathrm {1 \mskip-5.2mu l}}}}

\newcommand{\I}{\mathcal I}

\begin{document}
\title{An Indicator Function Limit Theorem in\\ Dynamical Systems}
\author{Olivier Durieu\thanks{Laboratoire de Math\'ematiques et Physique Th\'eorique, UMR 6083 CNRS, Universit\'e Fran\c{c}ois Rabelais, Tours;
e-mail: olivier.durieu@lmpt.univ-tours.fr}
 \, \& \, Dalibor Voln\'y\thanks{Laboratoire de Math\'ematiques Rapha\"el Salem,
 UMR 6085 CNRS, Universit\'e de Rouen;
e-mail: dalibor.volny@univ-rouen.fr}}
\maketitle

\begin{abstract}
We show by a constructive proof that
in all aperiodic dynamical system, for all sequences $(a_n)_{n\in\N}\subset\R_+$ such that $a_n\nearrow\infty$ and $\frac{a_n}{n}\rightarrow 0$ as $n\rightarrow\infty$, there exists a set $A\in\A$ having the property that
the sequence of the distributions of $(\frac{1}{a_{n}}S_{n}(\ind_A-\mu(A)))_{n\in\N}$ is dense in the space of all probability measures on $\R$. This extends the result of \cite{DurVol08b} to the non-ergodic case.

\medskip


\medskip

\end{abstract}

\section{Introduction and result}

Let $(\Omega,\A,\mu)$ be a probability space where $\Omega$ is a Lebesgue space and let $T$ be an invertible measure preserving transformation from $\Omega$ to $\Omega$. We say that $(\Omega,\A,\mu,T)$ is a dynamical system.
Further, the dynamical system is aperiodic if 
$$
\mu\{x\in\Omega\,:\,\exists n\ge 1, T^{n}x=x\}=0.
$$
It is ergodic if for any $A\in\A$, $T^{-1}A=A$ implies $\mu(A)=0$ or $1$.

\medskip

For a random variable $X$ from $\Omega$ to $\R$, we denote by $S_n(X)$ the partial sums $\sum_{i=0}^{n-1}X\circ T^i$, $n\ge 1$.

\medskip

The present paper concerns the question of the limit behavior of partial sums in general aperiodic dynamical systems.
In 1987, Burton and Denker \cite{BurDen87} proved that in any aperiodic dynamical system, there exists a function in $L_0^2$ which verifies the central limit theorem.
In general, for functions in $\mathbb{L}^p$ spaces,
Voln\'y \cite{Vol90} proved that for any sequence $a_n\to \infty$, $\frac{a_n}{n}\to 0$,
there exists a dense $G_\delta$ part $G$ of $L_0^p$ such that for any $f\in G$
the sequence of distributions of
$\frac{1}{a_{n_k}}S_{n_k}(f)$ is dense in the set of all probability measures on $\R$, see also Liardet and Voln\'y \cite{LiaVol97}. This work is also related to the question of the rate of convergence in the ergodic
theorem (see del Junco and  Rosenblatt \cite{JunRos79}).

\medskip

In Durieu and Voln\'y \cite{DurVol08b}, a similar result is shown for the class of centered indicator functions $\ind_A-\mu(A)$, $A\in\A$
and for ergodic dynamical systems. The following theorem was obtained.

\begin{theorem}\label{M2}
 Let $(\Omega,\A,\mu,T)$ be an ergodic dynamical system on a Lebesgue probability space, $(a_n)_{n\in\N}\subset\R_+$ be an increasing sequence
such that $a_n\nearrow\infty$ and $\frac{a_n}{n}\rightarrow 0$ as $n\rightarrow\infty$.

There exists a dense (for the pseudo-metric of the measure of the symmetric difference) $G_\delta$ class of sets $A\in\A$ having the property that for every probability $\nu$ on $\R$, there exists a subsequence $(n_k)_{k\in\N}$
satisfying
\begin{equation*}
 \frac{1}{a_{n_k}}S_{n_k}(\ind_A-\mu(A))\;\xrightarrow[\;k\to\infty\;]{\mathcal D}\;\nu.
\end{equation*}
\end{theorem}

\medskip

Here, we answer the question of the existence of a similar result in the non-ergodic case.

\medskip

Assume now that $(\Omega,\A,\mu,T)$ is not ergodic.
Let $(a_n)_{n\in\N}\subset\R_+$ be an increasing sequence
satisfying $a_n\nearrow\infty$ and $\frac{a_n}{n}\rightarrow 0$ as $n\rightarrow\infty$.
Denote by $\I$ the $\sigma$-algebra of the invariant sets and $(\mu^x)_{x\in\chi}$ the ergodic components of the measure $\mu$.
If there exist a set $A\in\A$, a probability measure $\nu$ on $\R$ and a sequence $(n_k)_{k\in\N}$ such that
$$
\frac{1}{a_{n_k}}S_{n_k}(\ind_A-\mu(A))\;\xrightarrow[\;k\to\infty\;]{\mathcal D}\;\nu,
$$
then $\esp(\ind_A|\I)=\mu(A)$ almost surely.

\medskip

Indeed, if there exists $x\in\chi$ such that $\mu^x(A)-\mu(A)=c>0$, then by Birkhoff's Ergodic Theorem,
$$
\frac{1}{n}S_{n}(\ind_A-\mu(A))\;\xrightarrow[\;n\to\infty\;]{}\;c
$$
$\mu^x$-almost surely. Therefore
$$
\frac{1}{a_{n_k}}S_{n_k}(\ind_A-\mu(A))\;\xrightarrow[\;k\to\infty\;]{}\;+\infty
$$
and we have a contradiction.

\medskip

So, to find a set which satisfies the conclusion of Theorem \ref{M2}, we have to consider the
sets $A$ such that $\esp(\ind_A|\I)$ is almost surely constant. The class of such sets is not, in general, dense in $\A$.
So, in the non-ergodic case, we cannot expect the result of genericity.

\medskip

Nevertheless, in the non-ergodic case, one can show the existence of an arbitrarily small set $A\in\A$ such that the sequence of the distributions of $\left(\frac{1}{a_{n}}S_{n}(\ind_A-\mu(A))\right)_{n\in\N}$ is dense in the set of probability measures on $\R$.

\medskip

We prove the following result.

\begin{theorem}\label{pppp}
Let $(\Omega,\A,\mu,T)$ be an aperiodic dynamical system on a Lebesgue probability space and $(a_n)_{n\in\N}\subset\R_+$ be an increasing sequence  
such that $a_n\nearrow\infty$ and $\frac{a_n}{n}\rightarrow 0$ as $n\rightarrow\infty$.
For all $\varepsilon>0$, there exists a set $A\in\A$ with $\mu(A)<\varepsilon$ such that for every probability measure $\nu$ on $\R$, there exists a sequence $(n_k)_{k\in\N}$
such that
\begin{equation*}
 \frac{1}{a_{n_k}}S_{n_k}(\ind_A-\mu(A))\;\xrightarrow[\;k\to\infty\;]{\mathcal D}\;\nu.
\end{equation*}
\end{theorem}

The rest of the paper is devoted to the proof of Theorem \ref{pppp}. Note that the proof that we propose is constructive.

\section{Proof}
Let $(\Omega,\A,\mu,T)$ be an aperiodic dynamical system and $(a_n)_{n\in\N}\subset\R_+$ be an increasing sequence  
such that $a_n\nearrow\infty$ and $\frac{a_n}{n}\rightarrow 0$ as $n\rightarrow\infty$ which are fixed for all the sequel.

\subsection{An equivalent statement}

Let $\M$ be the set of all probability measures on $\R$ and
$\M_0$ be the set of all probability measures on $\R$ which have zero-mean.
We denote by $d$ the L\'evy metric on $\M$. For all $\mu$ and $\nu$ in $\M$ with distribution functions
$F$ and $G$,
$$
d(\mu,\nu)=\inf\{\varepsilon>0\,:\,G(t-\varepsilon)-\varepsilon\le F(t)\le G(t+\varepsilon)+\varepsilon, \forall t\in\R\}.
$$
The space $(\M,d)$ is a complete separable metric space and convergence with respect to $d$ is equivalent to weak convergence of distributions (see Dudley \cite{Dud02}, pages 394-395).
If $X:\Omega\longrightarrow\R$ is a random variable, we denote by $\cL_\Omega(X)$ the distribution of $X$ on $\R$.
Using the separability of the set $\M_0$ which is dense in $\M$, we can prove that the next theorem is equivalent to Theorem \ref{pppp}.

\begin{theorem}\label{ppp}
For every $\varepsilon>0$ and for every sequence $(\nu_k)_{k\in\N}$ in $\M_0$, there exist a set $A\in\A$, with $\mu(A)<\varepsilon$, and a sequence $(n_k)_{k\in\N}$
such that
\begin{equation*}
 d(\cL_\Omega(\frac{1}{a_{n_k}}S_{n_k}(\ind_A-\mu(A))),\nu_k)\;\xrightarrow[\;k\to\infty\;]{}\;0.
\end{equation*}
\end{theorem}

\begin{proposition}

Theorem \ref{ppp} and Theorem \ref{pppp} are equivalent.

\end{proposition}

\proof

Let $M$ be a countable and dense subset of $\M_0$. We can find a sequence $(\nu_k)_{k\in\N}$ such that
for all $\eta\in M$, there exists an infinite set $K_\eta\subset\N$ verifying that for all $k\in K_\eta$,
$\nu_k=\eta$.

\medskip

Let $A$ be the set and  $(n_k)_{k\in\N}$ be the sequence associated to the sequence $(\nu_k)_{k\in\N}$ as in Theorem \ref{ppp}.

\medskip

For each $\nu\in M$, there exists an increasing sequence $(k_j)_{j\in\N}$ such that $\nu_{k_j}=\nu$ for all $j\in\N$.
By Theorem \ref{ppp}, 
\begin{equation*}
 d(\cL_\Omega(\frac{1}{a_{n_{k_j}}}S_{n_{k_j}}(\ind_A-\mu(A))),\nu)\;\xrightarrow[\;j\to\infty\;]{}\;0.
\end{equation*}
By classical argument, Theorem \ref{pppp} follows.

\medskip

The fact that Theorem \ref{pppp} implies Theorem \ref{ppp} is clear.
\cqfd

\medskip

Now to prove Theorem \ref{ppp}, we will construct explicitly the set $A$.
To do that, we will use the four following lemmas.

\medskip

\subsection{Auxiliary results}

Let $\nu$ be a probability on $\R$. For $B\in\B(\R)$ with $\nu(B)>0$, $\nu_B$ denotes the probability on $\R$ defined by
$
\nu_B(A)=\nu(B)^{-1}\nu(A\cap B).
$
For $x\in\R$, $\nu_x$ denotes the probability on $\R$ defined by
$
\nu_x(B)=\nu(\{xb\,/\,b\in B\}).
$

\begin{lemma}\label{Lnot}
Some properties of the L\'evy metric:
\begin{enumerate}[{\rm (i)}]
 \item For each probability $\nu$ on $\R$, for all Borel sets $B$,
$
d(\nu_B,\nu)\le \nu(\R\setminus B).
$
\item For all probabilities $\nu$ and $\eta$ on $\R$, for all $x\ge 1$,
$
d(\nu_x,\eta_x)\le d(\nu,\eta).
$
\item For all probability $\nu$ on $\R$, for all measurable functions $f$ and $g$ from $\Omega$ to $\R$,
$$
d(\L_\Omega(f+g),\nu)\le (\L_\Omega(f),\nu)+d(\L_\Omega(g),\delta_0)
$$
where $\delta_0$ is the Dirac measure at $0$.
\item For all probability $\nu$ on $\R$, 
$d(\nu,\delta_0)\le A$ if and only if $\nu((-\infty,-A))\le A$ and $\nu((A,\infty))\le A$.
\end{enumerate}
\end{lemma}
The proof is left to the reader.

\medskip

\begin{lemma}\label{L}
For all probability $\nu$ on $\R$, for all $\varepsilon>0$, there exists $C_0\ge 1$ and $n_0\in\N$, for all
$C\ge C_0$ and $n\ge n_0$, there exists a probability $\eta$ on $\R$ with support $S\subset [-a_nC,a_nC]\cap\Z$ such that for all $i\in S$, $\eta(\{i\})\in\Q$,
$d(\eta_{a_n},\nu)\le \varepsilon$
and
$\esp(\eta):=\int x {\rm d}\eta(x)=0$.
\end{lemma}

\proof

This Lemma is a consequence of Lemma 3.3 in \cite{DurVol08b} (which has a constructive proof) and of the fact that
for all probability measure on $\Z$ with finite support, we can find a probability on $\Z$ with same support which is arbitrarily close to the first one (with respect to $d$) and takes values in $\Q$. 
\cqfd

\medskip

The following lemma is classical, we do not give a proof.
\begin{lemma}\label{function}
 Let $(\Omega,\A,\mu)$ be a Lebesgue probability space and $\nu$ be a probability on $\R$.
Then, there exists a measurable random variable $X:\Omega\longrightarrow\R$, such that
$\cL_{\Omega}(X)=\nu$.
\end{lemma}

\medskip

Recall that a set $F\in\A$ is the base of a Rokhlin tower of height $n$ if the sets $F, TF, \dots, T^{n-1}F$ are pairwise disjoint.

\begin{lemma}\label{rokhlin}
For all $n\ge 1$ and for all $\varepsilon>0$, there exists
a set $F\in\A$ such that $\{F,\dots,T^{n-1}F\}$ is a Rokhlin tower of measure greater than $1-\varepsilon$ and the sojourn time 
in the junk set $J=\Omega\setminus(\cup_{i=0}^{n-1}T^iF)$ is almost surely $1$, i.e. 
for a.e. $x\in J$, $Tx\in F$.
\end{lemma}

\proof

This can be view as a consequence of Alpern's theorem \cite{Alp79}, by constructing a Rokhlin castle with two towers of height $n$ and $n+1$ and the base of the second tower of measure less than $\varepsilon$.
\cqfd

\medskip

\subsection{Proof of Theorem \ref{ppp}}

Let the sequence $(\nu_k)_{k\ge 1}$ in $\M_0$ and the constant $\varepsilon\in(0,1)$ be fixed.
Let $(\varepsilon_k)_{k\ge 1}$ be a decreasing sequence of positive reals such that $\sum_{k\ge 1}\varepsilon_k<\varepsilon$ and $\sum_{k\ge 1}k\varepsilon_k<\infty$.

\medskip

Theorem \ref{ppp} is a consequence of the following proposition, which is proved in the next section.
\begin{proposition}\label{Ak}
There exist a sequence of pairwise disjoint sets $A_k\in\A$ and a sequence of integers $(n_k)_{k\ge 1}$
such that,
\begin{enumerate}[{\rm (i)}]
\item $\mu(A_1)\le \varepsilon_1$ and for all $k> 1$, $\mu(A_k)\le \frac{a_{n_{k-1}}}{n_{k-1}}\varepsilon_k$;
\item
for all $k\ge 1$, 
$$
d(\cL_\Omega(\frac{1}{a_{n_k}}S_{n_k}(\ind_{A_k}-\mu(A_k))),\nu_k)\le \varepsilon_k;
$$
\item for all $k\ge 1$ and for all $j> k$,
$$
d(\cL_\Omega(\frac{1}{a_{n_j}}S_{n_j}(\ind_{A_k}-\mu(A_k))),\delta_0)\le \varepsilon_j.
$$
\end{enumerate}
\end{proposition}
We admit the proposition for the end of the proof.

Set
$$
A=\bigcup_{k\ge 1}A_k.
$$

Then by (i),
\[
\mu(A)=\sum_{k\ge 1}\mu(A_k)\le\sum_{k\ge 1}\varepsilon_k\le\varepsilon.
\]
By Proposition \ref{Ak} and by (iii) and (iv) of Lemma \ref{Lnot}, for all $j\ge 1$,
\begin{eqnarray*}
d(\cL_\Omega(\frac{1}{a_{n_j}}S_{n_j}(\ind_{A}-\mu(A))),\nu_j)
&\leq& \sum_{i=1}^{j-1}d(\cL_\Omega(\frac{1}{a_{n_j}}S_{n_j}(\ind_{A_i}-\mu(A_i))),\delta_0)\\
&&+d(\cL_\Omega(\frac{1}{a_{n_j}}S_{n_j}(\ind_{A_j}-\mu(A_j))),\nu_j)\\
&&+\sum_{i=j+1}^{\infty}d(\cL_\Omega(\frac{1}{a_{n_j}}S_{n_j}(\ind_{A_i}-\mu(A_i))),\delta_0)\\
&\le& (j-1)\varepsilon_j + \varepsilon_j+ \sum_{i=j+1}^{\infty}\frac{n_j}{a_{n_j}}\mu(A_i)\\
&\le& \sum_{i\ge j}i\varepsilon_i
\end{eqnarray*}
which goes to $0$ when $j$ goes to $\infty$.
Thus Theorem \ref{ppp} is proved.

\subsection{Proof of Proposition \ref{Ak}}

We give an explicit construction of the sets $A_k$. We begin by the construction of the set $A_1$.

\paragraph*{Step 1: the set $A_1$}\

The goal is to find a set $A_1$ and an integer $n_1$ such that 
$$
d(\cL_\Omega(\frac{1}{a_{n_1}}S_{n_1}(\ind_{A_1}-\mu(A_1))),\nu_1)\le \varepsilon_1.
$$
But we also want that the set $A_1$ becomes negligible for the partial sums of length $n_k$, $k\ge 2$ (condition (iii)).
There are several steps. First, we will define a set $A_{1,1}$ which satisfies (ii) and (iii) for $j=2$. Then, we will modify this set,
step by step, to have (iii) for all $j>2$.

\subparagraph*{The set $A_{1,1}$}\label{SA11}\

We consider the probability $\nu_1$, the constant $\varepsilon_1$ and
we set $\alpha_1:=\frac{\varepsilon_1}{8}$. Applying Lemma \ref{L} to $\nu_1$ and $\alpha_1$,
we get two constants $C(\nu_1,\alpha_1)$ and $n(\nu_1,\alpha_1)$ and we choose
$C_1:= C(\nu_1,\alpha_1)$ and $n_1\ge n(\nu_1,\alpha_1)$ such that
\begin{equation}\label{equ1}
\frac{d_1}{n_1}\le \alpha_1\quad\mbox{ where }d_1:=\lfloor a_{n_1}C_1\rfloor+1.
\end{equation}
We get a corresponding centered probability $\eta_1$ (given by Lemma \ref{L}) with support in $\{-d_1+1,\dots,d_1-1\}$
such that $d(\eta_{1\,a_{n_1}},\nu_1)\le\alpha_1$.
Since for all $i$, $\eta_1(\{i\})\in\Q$, there exist $q_1\in\N$ and $q_1^{(i)}\in\N$, $i=1\dots,2d_1$, such that
$$
\eta_1(\{i-d_1\})=\frac{q_1^{(i)}}{q_1},\,\mbox{ for all }i=1,\dots,2d_1.
$$

\bigskip

Now, we consider the probability $\nu_2$ and $\varepsilon_2$.
We define $\alpha_2:=\frac{a_{n_{1}}}{2n_{1}}\varepsilon_2$.
Applying Lemma \ref{L} to $\nu_2$ and $\alpha_2$,
we get two constants $C(\nu_2,\alpha_2)$ and $n(\nu_2,\alpha_2)$.
Set $C_2:= \max\{C(\nu_2,\alpha_2),C_{1}\}$ and let $n_2\ge n(\nu_2,\alpha_2)$ be a multiple of $q_1n_1$ such that
\begin{equation}\label{equ2'2}
\frac{q_{1}n_1}{a_{n_2}}\le \alpha_2.
\end{equation}

By Lemma \ref{rokhlin}, we can consider a set $F_1\in\A$ such that $\{F_1,TF_1,\dots,T^{n_2-1}F_1\}$ is a Rokhlin tower of height $n_2$,
with the sojourn time in the junk set almost surely equal to $1$ and the measure of the junk set smaller than $\gamma_1:=\min\{\frac{a_{n_2}}{n_2}\alpha_2,\alpha_1\}$.
Write
$F_1^l:=T^{ln_1}F_1$, $l=0,\dots,p_1-1$. We thus have $p_1:=\frac{n_{2}}{n_1}$ towers $\{F_1^l,\dots,T^{ln_1-1}F_1^l\}$ of height $n_1$. Notice that by definition of $n_2$, $p_1$ is a multiple of $q_1$.

\medskip

By Lemma \ref{function}, let $h_1$ be a measurable function from $F_1$ to $\Z$ such that $\cL_{F_1}(h_1)=\eta_1$ and denote by $g_1$ the positive function equal to $h_1+d_1$.
Let 
$$
A_{F_1,i}:=g_1^{-1}(\{i\}),\; i=1,\dots,2d_1.
$$
Now, for all $i=1,\dots, 2d_1$, let $\{A_{F_1,i,1},\dots,A_{F_1,i,q_1^{(i)}}\}$ be a partition of the set $A_{F_1,i}$
into sets of measure $\frac{1}{q_1}\mu(F_1)$. We thus have a partition of $F_1$ into
$$
\{A_{F_1,1,1},\dots,A_{F_1,1,q^{(1)}_1},A_{F_1,2,1},\dots,A_{F_1,2,q_1^{(2)}},\dots\dots\dots,A_{F_1,2d_1,1},\dots,A_{F_1,2d_1,q_1^{(2d_1)}}\}.
$$
By induction, we define partitions of $F_1^l$ for $l=1,\dots,p_1-1$ by setting
$$
A_{F_1^l,i,j}:=\left\{
\begin{array}{ll}
 T^{n_1}A_{F_1^{l-1},i,j+1} & \mbox{ if }1\le j\le q_1^{(i)}-1\\
 T^{n_1}A_{F_1^{l-1},i+1,1} & \mbox{ if }j=q_1^{(i)} \mbox{ and } i<2d_1\\
 T^{n_1}A_{F_1^{l-1},1,1} & \mbox{ if }j=q_1^{(2d_1)} \mbox{ and } i=2d_1
\end{array}
\right..
$$
For all $l=0,\dots,p_1-1$ and for all $i=1,\dots,2d_1$, we set
$$
A_{F_1^l,i}:=\bigcup_{j=1}^{q_1^{(i)}}A_{F_1^l,i,j}
$$
and
$$
A_{F_1^l}:=\bigcup_{i=1}^{2d_1}\bigcup_{k=0}^{i-1}T^kA_{F_1^l,i}.
$$
Remark that for any $l\in\{0,\dots,p_1-1\}$, for any $x\in F_1^l$  and for any $i\in\{1,\dots,2d_1\}$,
$$
x\in A_{F_1^l,i} \mbox{ if and only if } S_{n_1}(\ind_{A_{F_1^l}})(x)=i.
$$

\medskip

Now, we define the set $A_{1,1}$ as follows
$$
A_{1,1}:=\bigcup_{l=0}^{p_1-1}A_{F_1^l}.
$$
Remark that for any $x\in F_1$, $S_{n_1q_1}(\ind_{A_{1,1}})(x)=d_1q_1$ and since $p_1$ is a multiple of $q_1$,
$S_{n_2}(\ind_{A_{1,1}})(x)=d_1p_1$.  Moreover,
for each $k\in\{1,\dots,\frac{n_2}{n_1q_1}\}$, for any $x\in T^{kn_1q_1}F_1$, we also have $S_{n_1q_1}(\ind_{A_{1,1}})(x)=d_1q_1$.

\begin{lemma}\label{A11}
$\,$
\begin{enumerate}[{\rm (i)}]
\item $\mu(A_{1,1})\le\alpha_1$;
\item 
$\displaystyle
A_{1,1}\subset \bigcup_{j=0}^{p_1-1}\bigcup_{i=0}^{2d_1-1}T^{i+jn_1}F_1;
$
\item 
$\displaystyle
d(\cL_\Omega(\frac{1}{a_{n_1}}S_{n_1}(\ind_{A_{1,1}}-\mu(A_{1,1}))),\nu_1)\le \varepsilon_1;
$
\item
$\displaystyle
d(\cL_\Omega(\frac{1}{a_{n_2}}S_{n_2}(\ind_{A_{1,1}}-\mu(A_{1,1}))),\delta_0)\le \alpha_2.
$
\end{enumerate}
\end{lemma}

\proof

For each $l=0,\dots,p_1-1$, $\mu(A_{F_1^l})=\esp_{F_1}(g_1)\mu(F_1)=d_1\mu(F_1)$. Therefore, by (\ref{equ1}),
$$
\mu(A_{1,1})=\sum_{l=0}^{p_1-1}\mu(A_{F_1^l})=p_1d_1\mu(F_1)\le \frac{p_1}{n_2}d_1=\frac{d_1}{n_1}\le\alpha_1.
$$

\medskip

By construction, (ii) is clear.

\medskip

Let $\Omega_{1}:=\bigcup_{l=0}^{p_1-1}\bigcup_{i=0}^{n_1-2d_1-1}T^{-i}F_1^l$.
We have
$$
\cL_{\Omega_1}(S_{n_1}(\ind_{A_{1,1}}))=\cL_{F_1}(g_1)
$$
and since $\esp_{F_1}(g_1)=d_1$, by centering,
$$
\cL_{\Omega_1}(S_{n_1}(\ind_{A_{1,1}}-\mu(A_{1,1})))=\cL_{F_1}(h_1).
$$
Now, since $\gamma_1\le\alpha_1$,
$$
\mu(\Omega_1)= p_1(n_1-2d_1)\mu(F_1)\ge (n_2-2p_1d_1)\frac{(1-\gamma_1)}{n_2}=1-\gamma_1-\frac{2d_1}{n_1}\ge 1-3\alpha_1.
$$
Thus, by Lemma \ref{Lnot} (i),
$$
d(\cL_{\Omega}(S_{n_1}(\ind_{A_{1,1}}-\mu(A_{1,1}))),\cL_{F_1}(h_1))\le 3\alpha_1.
$$
and by Lemma \ref{Lnot} (ii),
$$
d(\cL_{\Omega}(\frac{1}{a_{n_1}}S_{n_1}(\ind_{A_{1,1}}-\mu(A_{1,1}))),\cL_{F_1}(\frac{h_1}{a_{n_1}}))\le 3\alpha_1.
$$
We infer, by triangular inequality, that
$$
d(\cL_{\Omega}(\frac{1}{a_{n_1}}S_{n_1}(\ind_{A_{1,1}}-\mu(A_{1,1}))),\nu_1)\le \varepsilon_1
$$
and (iii) is proved.

\medskip

Recall that $n_2$ is a multiple of $n_1q_1$ and, by definition of $A_{1,1}$, $S_{n_1q_1}(\ind_{A_{1,1}})(x)=d_1q_1$ whenever
$x$ belongs to one of the $T^{kn_1q_1}F_1$ for $k=0,\dots,\frac{n_2}{n_1q_1}$.
Since the sojourn time in the junk set is $1$, we infer that for any $x\in\Omega$,
$$
(p_1-q_1)d_1\le S_{n_2}(\ind_{A_{1,1}})(x)\le (p_1+q_1)d_1.
$$
Using $\mu(A_{1,1})=p_1d_1\mu(F_1)$, we get
$$
|S_{n_2}(\ind_{A_{1,1}}-\mu(A_{1,1}))|\le p_1d_1|1-n_2\mu(F_1)|+q_1d_1\le p_1d_1\gamma_1+q_1d_1.
$$
Thus, since $\gamma_1\le \frac{a_{n_2}}{n_2}\alpha_2$ and by (\ref{equ2'2}), we have
$$
\frac{1}{a_{n_2}}|S_{n_2}(\ind_{A_{1,1}}-\mu(A_{1,1}))|\le \frac{d_1}{n_1}\alpha_2+\frac{d_1}{n_1}\frac{q_1n_1}{a_{n_2}}\le
2\alpha_1\alpha_2\le\alpha_2
$$
and (iv) follows from application of Lemma \ref{Lnot} (iv).
\cqfd

\medskip

Of course, the set $A_{1,1}$ is not defined well enough to be negligible for higher partial sums. So, we need to modify a small part of $A_{1,1}$. Thus we introduce a sequense of sets $A_{1,k}$, $k\ge 2$. The set $A_{1,1}$ can be considered as a first version of the set $A_1$ and the $A_{1,k}$, $k\ge 2$ are the adjustments.

\subparagraph*{The sets $A_{1,k}$, $k\ge 2$}\ 

We shall give here the general algorithm to deduce the set $A_{1,k}$ from $A_{1,k-1}$. To do that, we need first to define
the entire sequence $(n_k)_{k\ge1}$.

\medskip

By induction, we define the sequences $(\alpha_k)_{k\ge 2}$, $(C_k)_{k\ge 2}$, $(n_k)_{k\ge 2}$, $(q_k)_{k\ge 2}$ as follows.
We consider the probability $\nu_k$ and $\varepsilon_k$.
We define $\alpha_k:=\frac{a_{n_{k-1}}}{2n_{k-1}}\varepsilon_k$.
Applying Lemma \ref{L} to $\nu_k$ and $\alpha_k$,
we get two constants $C(\nu_k,\alpha_k)$ and $n(\nu_k,\alpha_k)$.
Set $C_k:= \max\{C(\nu_k,\alpha_k),C_{k-1}\}$ and 
let $n_k\ge n(\nu_k,\alpha_k)$ be a multiple of $q_{k-1}n_{k-1}$ such that
\begin{equation}\label{equ1'}
\frac{d_k}{n_k}\le \alpha_k\quad\mbox{ where }d_k:=\lfloor a_{n_k}C_k\rfloor+1.
\end{equation}
and
\begin{equation}\label{equ2'}
\frac{q_{k-1}n_{k-1}}{a_{n_k}}\le \alpha_k.
\end{equation}
By Lemma \ref{L}, we get a corresponding centered probability $\eta_k$ on $\{-d_k+1,\dots,d_k-1\}$
such that $d(\eta_{k\,a_{n_k}},\nu_k)\le\alpha_k$.
There exist $q_k\in\N$ and $q_k^{(i)}\in\N$, $i=1\dots,2d_k$, such that
$\eta_k(\{i-d_k\})=\frac{q_k^{(i)}}{q_k}$.

\medskip

For $k\ge 1$, we also set $p_k:=\frac{n_{k+1}}{n_k}\in\N$ and $\beta_k:=\alpha_k-\alpha_{k+1}$. Thus, for all $k\ge1$,
\begin{equation}\label{equ3}
 \sum_{j\ge k}\beta_j\le\alpha_k.
\end{equation}
We define the sequence $(\gamma_k)_{k\ge 1}$ by
\begin{equation}\label{equ4}
\gamma_k:= \min\left\{ \frac{\beta_{k+1}}{2p_{k+1}}, \frac{a_{n_{k+1}}}{n_{k+1}}\alpha_{k+1}\right\}.
\end{equation}

Further, for all $k\ge 1$, by application of Lemma \ref{rokhlin}, we obtain a set $F_k\in\A$ such that
$\{F_k,TF_k,\dots,T^{n_{k+1}-1}F_k\}$ is a Rokhlin tower of height $n_{k+1}$ and the
junk set $J_k:=\Omega\setminus\bigcup_{i=0}^{n_{k+1}-1}T^iF_k$ is a set with sojourn time 1 and
$\mu(J_k)\le\gamma_k$.

\medskip

Remark that $\alpha_2$, $n_2$ and $\gamma_2$ have been previously defined but they respect this new definition.

\medskip

We also introduce the sequence of sets $F_k'$ defined by induction in the following way:
$F_1':=F_1$ and 
$$
F_k':=\bigcup_{x\in F_k}T^{n(x)}x.
$$
where for all $x$ in $F_k$, $n(x):=\inf\{n\ge 0\,/\,T^nx\in F_{k-1}'\}$ is the time of the first visit in $F_{k-1}'$.

\bigskip

\begin{lemma}\label{A1k}
There exists a sequence of measurable sets $(A_{1,k})_{k\ge 1}$ such that
\begin{enumerate}[{\rm (i)}]
\item
$\displaystyle
\mu(A_{1,k-1}\triangle A_{1,k})\le \beta_k;
$
\item for all $x\in\Omega$, $S_{n_1}(\ind_{A_{1,k}})(x)\le 2d_1$;
\item for all $x\in F_k'$, for all $i=0,\dots,p_k-1$, $S_{n_{k}}(\ind_{A_{1,k}}\circ T^{in_k})(x)=\frac{n_{k}}{n_1}d_1$;
\item
$\displaystyle
d(\cL_\Omega(\frac{1}{a_{n_{k+1}}}S_{n_{k+1}}(\ind_{A_{1,k}}-\mu(A_{1,k})),\delta_0)\leq\alpha_{k+1}.
$
 \end{enumerate}
\end{lemma}

\proof

We prove the lemma by induction. The set $A_{1,1}$ is already defined.

\medskip

Now, for a fixed $k$, we are going to explain how to deduce the set $A_{1,k}$ from $A_{1,k-1}$.

\medskip

For $x\in F_k'$ and $i=0,\dots,p_k-1$, let 
$$
\rho_i(x):=\sum_{j=in_k}^{(i+1)n_k-1}\ind_{A_{1,k-1}}\circ T^j(x)=S_{n_k}(\ind_{A_{1,k-1}}\circ T^{in_k}).
$$

By hypothesis, 
for all $x\in F_k'$, $\rho_0(x)=S_{n_k}(\ind_{A_{1,k-1}})(x)=p_{k-1}\frac{n_{k-1}}{n_1}d_1=\frac{n_k}{n_1}d_1$ 
but for $i>0$ it can be different. The differences appear when the orbit of the point $x$
meets the junk set $J_{k-1}$. Nevertheless, by definition of the Rokhlin tower (see Lemma \ref{rokhlin}), it can meet $J_{k-1}$ only one time in every $n_k$ consecutive iterates by $T$. 
So we have, $(p_{k-1}-1)\frac{n_{k-1}}{n_1}d_1 \le\rho_i(x)\le (p_{k-1}+1)\frac{n_{k-1}}{n_1}d_1$.
To summarize, for $x\in F_k'$ and $i=0,\dots,p_k-1$, 
$$
\rho_i(x)=\frac{n_k}{n_1}d_1+j
$$
with
$j\in\{-\frac{n_{k-1}}{n_1}d_1,\dots,\frac{n_{k-1}}{n_1}d_1\}$.

\medskip

We define a set $B_i(x)$ as follows. If $j\ge 0$, $B_i(x)=\emptyset$. If 
$j<0$, let $B_i(x)$ be a set composed by $|j|$ points from the set $\{T^{in_k}x,\dots,T^{(i+1)n_k-1}x\}\setminus A_{1,k-1}$,
in such a way that every $n_1$-consecutive points in  $\{T^{in_k}x,\dots,T^{(i+1)n_k-1}x\}$ meet $A_{1,k-1}\cup B_i(x)$ at most $2d_1$ times (it is possible because otherwise, $\rho_i(x)>\frac{n_{k}}{n_1}2d_1-|j|\ge\frac{n_k}{n_1}d_1$).

\medskip

We define a set $C_i(x)$ as follows. If $j\le 0$, $C_i(x)=\emptyset$. 
If $j>0$, let $C_i(x)$ be the set composed by the $j$ first points of $\{T^{in_k}x,\dots\}\cap
A_{1,k-1}$.

\medskip

Let 
$$
B:=\bigcup_{x\in F_k'}\bigcup_{i=0}^{p_k-1}B_i(x)\quad,\quad
C:=\bigcup_{x\in F_k'}\bigcup_{i=0}^{p_k-1}C_i(x)
$$
and
$$
A_{1,k}:=(A_{1,k-1}\setminus C)\cup B.
$$

\bigskip

Since the orbit of a point $x$ can only meet $J_{k-1}$ one time every $n_k$ and using (\ref{equ4}), we have
$$
\mu(A_{1,k-1}\triangle A_{1,k})\le 2p_k\mu(J_{k-1})\le 2p_k\gamma_{k-1}\le\beta_k.
$$

\medskip

Remark that (ii) and (iii) are guaranteed by construction of $A_{1,k}$.

\medskip

Further for all $x\in F_k'$, we have
$$
S_{n_{k+1}}(\ind_{A_{1,k}})(x)=p_k\frac{n_{k}}{n_1}d_1=\frac{n_{k+1}}{n_1}d_1.
$$
We deduce that $|\mu(A_{1,k})-\frac{n_{k+1}}{n_1}d_1\mu(F_k')|\le \mu(J_k)$ and
$$
\left(p_k-1\right)\frac{n_{k}}{n_1}d_1\le S_{n_{k+1}}(\ind_{A_{1,k}})\le \left(p_k+1\right)\frac{n_{k}}{n_1}d_1.
$$
Then,
$$
|S_{n_{k+1}}(\ind_{A_{1,k}}-\mu(A_{1,k}))|\le\frac{n_{k}}{n_1}d_1+(1+\frac{n_{k+1}}{n_1}d_1)\mu(J_k)
$$
and  by (\ref{equ4}) and (\ref{equ2'}),
$$
\frac{1}{a_{n_{k+1}}}|S_{n_{k+1}}(\ind_{A_{1,k}}-\mu(A_{1,k}))|\le
\frac{n_k}{a_{n_{k+1}}}\frac{d_1}{n_1}+(\frac{1}{n_{k+1}}+\frac{d_1}{n_1})\alpha_{k+1}
\le \alpha_{k+1}.
$$
By Lemma \ref{Lnot} (iv), we get (iv).
\cqfd

\medskip

\subparagraph*{The set $A_1$}\

We can now define the set $A_1\in\A$ as
$$
A_1:=\lim_{k\rightarrow\infty}A_{1,k}
$$
which is well defined because the sequence $(\mu(A_{1,k}\triangle A_{1,k+1}))_{k\ge  1}$ is summable. 
\begin{lemma}\label{A1}
$\;$
\begin{enumerate}[{\rm (i)}]
\item $\mu(A_1)\le 2\alpha_1$;
\item $S_{n_1}(\ind_{A_1})\le 2d_1$;
\item 
$\displaystyle
d(\cL_\Omega(\frac{1}{a_{n_1}}S_{n_1}(\ind_{A_{1}}-\mu(A_{1}))),\nu_1)\le 2\varepsilon_1;
$
\item For all $k\ge2$,
$\displaystyle
d(\cL_\Omega(\frac{1}{a_{n_k}}S_{n_k}(\ind_{A_{1}}-\mu(A_{1}))),\delta_0)\leq\varepsilon_k.
$
 \end{enumerate}
\end{lemma}

\proof

For all $k\ge1$, we have
\begin{equation}\label{difsym}
\mu(A_{1}\triangle A_{1,k})\le \sum_{j=k+1}^\infty\mu(A_{1,j-1}\triangle A_{1,j})\le \sum_{j=k+1}^\infty\beta_j\le\alpha_{k+1}
\end{equation}
and then $\mu(A_1)\le\mu(A_{1,1})+\mu(A_{1}\triangle A_{1,1})\le2\alpha_1$.

\medskip

Assertion (ii) comes from Lemma \ref{A1k} (ii).

\medskip

Further (\ref{difsym}) and Lemma \ref{Lnot} (iv) imply that for all $n$,
$$
d(\cL_\Omega(\frac{1}{a_{n}}(S_{n}(\ind_{A_{1,k}}-\mu(A_{1,k}))-S_{n}(\ind_{A_{1}}-\mu(A_{1})))),\delta_0)\le
\frac{n}{a_n}\alpha_{k+1}.
$$
Using Lemma \ref{Lnot} (iii), we can deduce (ii) from Lemma \ref{A11} (iii) and (iii) from Lemma \ref{A1k} (iv).
\cqfd

\medskip

\paragraph*{Step 2: The set $A_2$}\

\subparagraph*{The set $A_{2,2}$}\

We consider $F_2\in\A$ and we will almost repeat what we did to find the set $A_{1,1}$, working with $n_2, q_2, p_2, d_2$ instead
of $n_1, q_1, p_1, d_1$. The difference comes to the fact that we want $A_1\cap A_2=\emptyset$.
Recall that $\eta_2$ is the probability measure with support in $\Z$ given by Lemma \ref{L} applied to $\nu_2$ and $\alpha_2$ and with constants $C_2$ and $n_2$.
Let $h_2$ be a function from $F_2$ to $\Z$ given by Lemma \ref{function} such that $\cL_G(h_2)=\eta_2$ and call $g_2$ the positive function equal to $h_2+d_2$.
Let 
$$
A_{F_2,i}:=g_2^{-1}(\{i\}),\; i=1,\dots,2d_2
$$
and for all $i=1,\dots, 2d_2$, let $\{A_{F_2,i,1},\dots,A_{F_2,i,q_2^{(i)}}\}$ be a partition of the set $A_{F_2,i}$
into sets of measure $\frac{1}{q_2}\mu(F_2)$. We thus have a partition of $F_2$ into
$$
\{A_{F_2,1,1},\dots,A_{F_2,1,q^{(1)}_2},A_{F_2,2,1},\dots,A_{F_2,2,q_2^{(2)}},\dots\dots\dots,A_{F_2,2d_1,1},\dots,A_{F_2,2d_1,q_2^{(2d_1)}}\}.
$$
By induction, we deduce partitions of $F_2^l=T^{ln_2}F_2$ for $l=1,\dots,p_2-1$. We set
$$
A_{F_2^l,i,j}:=\left\{
\begin{array}{ll}
 T^{n_2}A_{F_2^{l-1},i,j+1} & \mbox{ if }1\le j\le q_2^{(i)}-1\\
 T^{n_2}A_{F_2^{l-1},i+1,1} & \mbox{ if }j=q_2^{(i)} \mbox{ and } i<2d_2\\
 T^{n_2}A_{F_2^{l-1},1,1} & \mbox{ if }j=q_2^{(2d_1)} \mbox{ and } i=2d_2
\end{array}
\right..
$$
For all $l=0,\dots,p_2-1$ and for all $i=1,\dots,2d_2$, we set
$$
A_{F_2^l,i}:=\bigcup_{j=1}^{q_2^{(i)}}A_{F_2^l,i,j}.
$$
Because we want disjointness, we cannot define $A_{F_2^l}$ as $\bigcup_{i=1}^{2d_2}\bigcup_{k=0}^{i-1}T^kA_{F_2^l,i}$.
So, for each $l\in\{0,\dots,p_2-1\}$, for each $x\in F_2^l$, if $x\in A_{F_2^l,i}$, we denote by $D_l(x)$ the set composed by the $i$ first elements of $\{x,Tx,\dots\}\setminus A_1$ and we set
$$
D_l:=\bigcup_{x\in F_2^l}D_l(x).
$$
Since $A_1$ contains at most $2d_1$ points in each orbit
of size $n_1$ and since $d_2\ge d_1$, 
$$
D_l\subset\bigcup_{i=0}^{4d_2-1}T^iF_2^l.
$$ 

We define $A_{2,2}$ as
$$
A_{2,2}:=\bigcup_{l=0}^{p_2-1}D_l.
$$
\begin{lemma}\label{A22}
$\,$
\begin{enumerate}[{\rm (i)}]
\item $\mu(A_{2,2})\le\alpha_2$;
\item  $S_{n_2}(\ind_{A_{2,2}})\le 2d_2$;
\item 
$\displaystyle
d(\cL_\Omega(\frac{1}{a_{n_2}}S_{n_2}(\ind_{A_{2,2}}-\mu(A_{2,2}))),\nu_2)\le \varepsilon_2;
$
\item
$\displaystyle
d(\cL_\Omega(\frac{1}{a_{n_3}}S_{n_3}(\ind_{A_{2,2}}-\mu(A_{2,2}))),\delta_0)\le \alpha_3.
$
\end{enumerate}
\end{lemma}

\proof

The proof follows the one of Lemma \ref{A11} and is left to the reader.
\cqfd

\medskip

\subparagraph*{The sets $A_{2,k}$, $k\ge 3$}\

Now we define a sequence $A_{2,k},k\ge 3$
using the same techniques as before and preserving the fact that for all $k\ge 1$,
$A_{1,k}$ contains at most $2d_2$ points in each orbit of length $n_2$.

Then the $A_{2,k}$ satisfy
$$
\mu(A_{2,k-1}\triangle A_{2,k})\le \beta_{k}
$$
and
$$
d(\cL_\Omega(\frac{1}{a_{n_{k+1}}}S_{n_{k+1}}(\ind_{A_{2,k}}-\mu(A_{2,k}))),\delta_0)\le \alpha_{k+1}.
$$

\subparagraph*{The set $A_2$}\

The set
$$
A_2:=\lim_{k\rightarrow\infty}A_{2,k}
$$
is well defined, disjoint of $A_1$ and satisfies the following lemma.

\begin{lemma}\label{A2}
$\,$
\begin{enumerate}[{\rm (i)}]
\item $\mu(A_2)\le 2\alpha_2$;
\item  $S_{n_2}(\ind_{A_2})\le 2d_2$;
\item
$\displaystyle
d(\cL_\Omega(\frac{1}{a_{n_2}}S_{n_2}(\ind_{A_2}-\mu(A_2))),\nu_2)\le 2\varepsilon_2;
$
\item For all $k\ge 3$,
$\displaystyle
d(\cL_\Omega(\frac{1}{a_{n_k}}S_{n_k}(\ind_{A_2}-\mu(A_2))),\delta_0)\le \varepsilon_k.
$
\end{enumerate}
\end{lemma}

\proof

The proof follows the one of Lemma \ref{A1}.
\cqfd

\medskip

\paragraph*{Step k: the set $A_k$, $k\ge3$}\

It is now clear that, by induction, we can find sets $A_k$ disjoint of the sets $A_i$, $i<k$, satisfying Proposition \ref{Ak}.
\cqfd

\bibliographystyle{plain}

\begin{thebibliography}{1}

\bibitem{Alp79}
Steve Alpern.
\newblock Generic properties of measure preserving homeomorphisms.
\newblock In {\em Ergodic theory (Proc. Conf., Math. Forschungsinst.,
  Oberwolfach, 1978)}, volume 729 of {\em Lecture Notes in Math.}, pages
  16--27. Springer, Berlin, 1979.

\bibitem{BurDen87}
Robert Burton and Manfred Denker.
\newblock On the central limit theorem for dynamical systems.
\newblock {\em Trans. Amer. Math. Soc.}, 302(2):715--726, 1987.

\bibitem{JunRos79}
Andr{\'e}s del Junco and Joseph Rosenblatt.
\newblock Counterexamples in ergodic theory and number theory.
\newblock {\em Math. Ann.}, 245(3):185--197, 1979.

\bibitem{Dud02}
R.~M. Dudley.
\newblock {\em Real analysis and probability}, volume~74 of {\em Cambridge
  Studies in Advanced Mathematics}.
\newblock Cambridge University Press, Cambridge, 2002.
\newblock Revised reprint of the 1989 original.

\bibitem{DurVol08b}
Olivier Durieu and Dalibor Voln\'y.
\newblock On sums of indicator functions in dynamical systems.
\newblock {\em to appear in \textit{Ergodic Theory and Dynamical Systems}.}

\bibitem{LiaVol97}
Pierre Liardet and Dalibor Voln{\'y}.
\newblock Sums of continuous and differentiable functions in dynamical systems.
\newblock {\em Israel J. Math.}, 98:29--60, 1997.

\bibitem{Vol90}
Dalibor Voln{\'y}.
\newblock On limit theorems and category for dynamical systems.
\newblock {\em Yokohama Math. J.}, 38(1):29--35, 1990.

\end{thebibliography}

\end{document}